# Estimation of non-symmetric and unbounded region of attraction using shifted shape function and R-composition

Dongyang Li[a,b,*,1], Dmitry Ignatyev[b,2], Antonios Tsourdos[b,2], Zhongyuan Wang[a,1]

[a] *School of Energy and Power Engineering, Nanjing University of Science and Technology, Nanjing 210094, China*
[b] *School of Aerospace, Transport and Manufacturing, Cranfield University, Cranfield MK43 OAL, UK*

## Abstract

A general numerical method using sum of squares programming is proposed to address the problem of estimating the region of attraction (ROA) of an asymptotically stable equilibrium point of a nonlinear polynomial system. The method is based on Lyapunov theory, and a shape function is defined to enlarge the provable subset of a local Lyapunov function. In contrast with existing methods with a shape function centered at the equilibrium point, the proposed method utilizes a shifted shape function (SSF) with its center shifted iteratively towards the boundary of the newly obtained invariant subset to improve ROA estimation. A set of shifting centers with corresponding SSFs is generated to produce proven subsets of the exact ROA and then a composition method, namely R-composition, is employed to express these independent sets in a compact form by just a single but richer-shaped level set. The proposed method denoted as RcomSSF brings a significant improvement for general ROA estimation problems, especially for non-symmetric or unbounded ROA, while keeping the computational burden at a reasonable level. Its effectiveness and advantages are demonstrated by several benchmark examples from literature.

*Corresponding author
 E-mail address: dongylee@njust.edu.cn
 [1] Postal address: School of Energy and Power Engineering, Nanjing University of Science and Technology, Xuanwu District, Nanjing 210094, China.
  [2] Postal address: School of Aerospace, Transport and Manufacturing, Cranfield University, Cranfield MK43 OAL, UK.

# 1 Introduction

Many real-world systems are governed by systems of nonlinear equations [1]. Performing stability analysis analytically for a nonlinear system in general form is hard, if not impossible. That is why there is a tremendous interest in numerical methods allowing stability analysis. As opposed to Linear Time Invariant (LTI) systems, whose stability can be analyzed globally, most of the nonlinear systems achieve stability only in a specific region around some equilibrium point (EP) [1]. For practical applications, it is of great interest to find this specific invariant set [2,3]. This invariant set, also the so-called region of attraction (ROA) of the relevant EP, is an important metric for system stability and robustness, showing how much the initial states can be disturbed away from the expected steady states. Besides, an actual dynamical system can have multiple stable EPs (or limit cycles) [2]. Therefore, one has to ensure the operative range of the system is contained in the ROA of the expected EP. However, finding the exact ROA, both numerically and analytically, is hard [3]. Extensive efforts have been made in the last few years to obtain a closer approximation to the exact ROA in different applications, such as aerospace, robotics, medical, chemical process, traffics, and biological systems [2–24].

Among all methods for ROA estimation, analytical or computational, those based on Lyapunov function (LF) are the most popular ones. They provide sufficient conditions for the stability of equilibria. However, constructing LFs for a general system is not an easy task and has been studied extensively [25]. Many computational construction methods have been developed and various forms of relaxation are used. One of them is the sum of squares (SOS) technique. The condition for a multivariable polynomial being non-negative is relaxed by a SOS polynomial [20]. For example, $s(x)$ is a SOS polynomial if there exist polynomials $\{f_i(x)\}_{i=1}^{N} \in \mathbf{R}(x)$ such that $s(x) = \sum_{i=1}^{N} f_i^2(x)$, where $\mathbf{R}[x]$ represents the set of polynomials in $x \in \mathbb{R}^n$ with real coefficients and $N$ is a positive integer [20]. The relaxed SOS constraint is tighter but has the advantage of making tractable not only the construction of LF but also the computation of its invariant sublevel set, which is then useful for ROA estimation. Readers are referred to [12,20,21,26] for a more detailed discussion of SOS techniques.

SOS problems containing SOS constraints are treated with SOS programming. Freely available MATLAB toolboxes such as SOSTOOLs [27], SOSOPTs [28] transform SOS constraints into semidefinite problems (SDPs) and then are solved by SDP solvers such as SeDuMi [29]. For example, in SOS programming, verifying a polynomial $s(x)$ being a SOS polynomial is equivalent to checking the existence of a positive semidefinite matrix $Q$, such that

$$s(x) = Z^T(x)QZ(x) \in \Sigma_n$$

where $Z(x)$ is some properly chosen vector of monomials in $x \in \mathbb{R}^n$ and $\Sigma_n$ denotes the set of SOS polynomials in $x \in \mathbb{R}^n$. It is a basic feasibility problem in SOS programming and another common class is the optimization problem for a linear objective function [26]. These problems are formulated by linear matrix inequalities (LMIs) and solved by LMI solvers. However, the formulation of some practical problems might result in bilinear matrix inequalities (BMIs) and cannot be solved easily by LMI solvers, which happens in finding the largest possible invariant subset of ROA [2,18]. In the ROA problem, decision variables being optimized are coupled with auxiliary SOS multipliers. One solution is using bilinear solvers (such as PENMBI in YALMIP [30]). Although they allow direct treatment of bilinear problems, they are less developed than the linear ones and the convergence to a global optimum cannot be guaranteed [18,21,22,31,32]. The other way is to convert it into a two-way iterative search between LFs and SOS multipliers by using the additional structure of the ROA estimation problem than a general bilinear problem.

Then it becomes affine in SOS multipliers if the LF is fixed, and vice versa [27,33]. The two-way search can then be solved conveniently by the widely used linear solvers.

The two-way iteration procedure has been used for many applications and studied for improvements in literature. For example, an interior expanding algorithm [34] is proposed to enlarge the ROA estimation using a positive definite polynomial (or 'shape function') and the composed algorithm is denoted as V-s iteration algorithm [23,35,36]. Elements that affect the performance of V-s iteration, e.g. choices of initial LF $V_0$, the degree of the new LF $V$, and the shape function $p(x)$, have been considered to achieve a larger estimation. A higher-degree LF is searched rather than the usual quadratic ones in order for a richer level set [31]. Methods are proposed to address the dramatic increase in computation cost with an increase of system dimension and/or polynomial degree. Thus, a composite LF by lower degree LFs, such as the pointwise maximum and pointwise minimum of LFs, are employed [18,21,37]. Later, rather than via only one Lyapunov estimate, a family of parameter-dependent LFs is used for an improved estimation [38]. A systematic way (denoted as R-composition [39]) is applied for a richer-shaped estimation [40] by composing LFs through the use of R-functions. Optimization of the initial LF $V_0$ is also important. The method in [5] proposes to optimize $V_0$ to avoid numerical infeasibility at the initial step when the level set is lower than solver tolerance. Topcu et al. [18] propose to use information from simulations to generate better LF candidates to improve the performance of bilinear SDP solver. Regarding the shape function, a quadratic form is customized to align better with the shape of the exact ROA [2] or reflects the relative importance of states in some practical problems [41–43]. However, a proper shape function is difficult to decide without prior knowledge of the system, and no systematic approach is ever proposed [5,13,18] but a quadratic form of $p(x) = x^T N x$ is widely used with a general assumption of EP at the origin. The shape matrix

$N \in \mathbb{R}^{n \times n}$ ($n$ is the system dimension) is a positive definite matrix, which is problem-dependent and commonly used to scale the state space domain. In most applications, the shape function is set fixed throughout iteration [44]. While in [2], an adaptive shape function is proposed to update by the quadratic parts of the newly found LF each iteration. A better estimation is obtained for some examples, especially those with symmetric or simple-shaped ROA. However, it faces the same problem of convergence after certain iterations as the conventional V-s algorithm does and sometimes underperforms the fixed shape function algorithm.

In other ROA estimation methods, e.g. [5], an additional constraint was introduced to improve estimation but can be over-restrictive when a higher-degree LF is searched. Tan [31] proposed to use a series of consecutive rotations of an elliptical shape function for enlargement. A similar approach by varying the matrix $P$ is proposed in [22]. In [45], sufficient conditions are provided to guarantee the sublevel sets of polynomial LFs can inner approximate the exact ROA up to any desired accuracy, but only applicable to the bounded ROA. Its practical application using SOS programming on the Van de Pol system only gives a similar result as in [2]. A simple algorithm to handle the bilinear problem for ROA estimation is explained in [15] but it is essentially similar with the two-way iterative search in the widely used V-s iteration algorithm.

In addition to the modifications proposed above, V-s iteration algorithm uses an increased number of iterations as a means of enlarging ROA estimation. However, it is not a universal remedy, for example, under circumstances of converged optimization and numerical infeasibility, especially for a system with unbounded or irregular ROA. Another feature is that the geometric center of the shape function, either fixed or adaptive, is located at the origin, which limits the estimation, especially for non-symmetric or unbounded ROA.

The main contribution of this work is the proposal to use shifted shape functions with shifting centers to enlarge ROA estimation even for the non-symmetric or unbounded ROA. The form of the shifted shape function proposed in this paper has never been explored, to the best knowledge of the authors, though it can be absorbed in the general construction of SOS polynomials. The shifted shape functions are constructed by shifting centers chosen iteratively from the newly obtained proven level set of the true ROA. It implies the advantage of making full use of the results already obtained. The proven level set can be conveniently obtained by the V-s iteration algorithm. In addition, the V-s algorithm also enables the shifted shape function to satisfy the constraints. A set of shifting centers with corresponding shifted shape functions is generated to produce several proven subsets of the ROA. A union of these sets by R-composition gives a compact and richer-shaped result without an increase in the LF degree thus alleviating the computational burden. Combined, the proposed approach is denoted as RcomSSF in the remainder of the paper. It helps to push forward the extension when additional iterations do not give further ROA expansion or when numerical infeasibility is encountered. In addition, it can take full advantage of the existing improvements for V-s iteration algorithm, for instance, modifications on shape matrix and initial LF.

The paper is organized as follows. Firstly, the ROA estimation problem is described with fundamental theories in Sec. 2. Then by using SOS techniques, the computational V-s iteration algorithm is presented. Sec. 3 goes into detail about the shape function to provide the premise for the proposed method. Then comes the proposed RcomSSF in Sec. 4 with the construction and implementation details of the shifted shape function and R-composition. Sec. 5 illustrates several numerical examples from literature. By comparing with existing ROA estimation algorithms, the effectiveness and advantages of the proposed method are demonstrated. Finally, the paper concludes in Sec. 6.

## 2 Region of attraction estimation

### 2.1 Problem formulation

Consider an autonomous nonlinear polynomial system

$$\dot{x} = f(x), \quad x(0) = x_0, \tag{1}$$

where $x \in \mathbb{R}^n$ is the state vector and $f(x)$ is a $n \times 1$ polynomial vector field. Without loss of generality, we assume that the origin is an asymptotically stable EP such that $f(0) = 0$. The ROA, a set of initial conditions whose trajectories will not go beyond this set and always converge back to the origin, can be defined as

$$\Omega := \left\{ x_0 \in \mathbb{R}^n : \text{If } x(0) = x_0 \text{ then } \lim_{t \to \infty} x(t) = 0 \right\}. \tag{2}$$

Determining the exact ROA of the nonlinear system is hard if not impossible, which brings the problem of finding the approximation of the ROA and it has been studied extensively in literature [2,4–12,14,18–24,46]. The estimation of ROA in this paper is based on the following proved Lemma 1, following from the direct Lyapunov theorem that specifies a sublevel set of a LF as an inner approximation for the true ROA of an asymptotically stable EP [37].

**Lemma 1** If there exist a continuously differentiable scalar function $V(x): \mathbb{R}^n \to \mathbb{R}$ and a positive scalar $\gamma \in \mathbb{R}^+$, such that

$$V(x) > 0 \quad \forall x \neq 0 \text{ and } V(0) = 0, \tag{3}$$

$$\Omega_\gamma := \{x : V(x) \leq \gamma\} \text{ is bounded}, \tag{4}$$

$$\Omega_\gamma \subseteq \{x : (\partial V / \partial x) f < 0\} \cup \{0\}, \tag{5}$$

then the origin is asymptotically stable and $\Omega_\gamma$ is a subset of the ROA. The level $\gamma$ can be optimized to get the largest possible estimation. $V(x)$ that satisfies Lemma 1 is called a strict LF. When $\gamma$ in Eqn. (4) is unbounded, the system is globally asymptotically stable.

To enlarge the estimation $\Omega_\gamma$, the interior expanding algorithm [2,34,44] is introduced, where a scalar polynomial function $p(x) \in \mathbf{R}[x]$ ($\mathbf{R}[x]$ is the set of polynomials with real coefficients in $x \in \mathbb{R}^n$) and another positive scalar $\beta \in \mathbb{R}^+$ are defined. It is ensured that

$$\varepsilon_\beta := \{x \in \mathbb{R}^n : p(x) \leq \beta\} \text{ is bounded,} \tag{6}$$

$$\varepsilon_\beta \subseteq \Omega_\gamma. \tag{7}$$

$p(x)$, denoted as shape function, is a positive definite and convex polynomial. The positive scalar $\beta$ is maximized while imposing constraints (3)-(5).

Polynomials constraints above can be relaxed into SOS constraints by SOS techniques using the connection between nonnegativity and sum of squares [12,20] and then solved by SOS programming. For the set containment constraints in Eqn. (5) and (7), a well-known generalized S-procedure [18,20] will be employed.

**Lemma 2** (Generalized S-procedure [18]). Given polynomials $g_0(x), g_1(x)...g_m(x) \in \mathbf{R}[x]$ and polynomials $s_1(x)...s_m(x) \in \Sigma_n$, if

$$g_0(x) - \sum_{i=1}^m s_i(x)g_i(x) \geq 0, \tag{8}$$

then

$$\{x \mid g_1(x), g_2(x)...g_m(x) \geq 0\} \subseteq \{x \mid g_0(x) \geq 0\}. \tag{9}$$

Therefore, the set containment constraint in Eqn. (5) can be formulated as

$$-[(\partial V/\partial x)f + l_2] - (\gamma - V)s_2 \geq 0, \tag{10}$$

where $\Sigma_n$ denotes the set of SOS polynomials in $x \in \mathbb{R}^n$ respectively. And constraint (7) can be reformulated in the same manner. Therefore the above ROA estimation problem (3)-(7) can be formulated as an optimization problem

$$\begin{aligned} & \max_{s_1, s_2 \in \Sigma_n} \beta \\ \text{Subject to:} \quad & V - l_1 \in \Sigma_n \\ & -[(\partial V/\partial x)f + l_2] - (\gamma - V)s_2 \in \Sigma_n \\ & (\gamma - V) - (\beta - p)s_1 \in \Sigma_n \end{aligned}, \tag{11}$$

where $l_i(x)(i=1,2)$ is a small positive polynomial (typically $\varepsilon x^T x$ with some small $\varepsilon \in \mathbb{R}^+$) to guarantee the derivative of $V$ strictly negative; $s_i(x)(i=1,2)$ is an auxiliary SOS multiplier with a proper degree. The resulting optimization problem (11) is a bilinear problem with $\beta$, $\gamma$ and decision variables in $V$ coupled with that in $s_i$.

## 2.2 V-s **iteration algorithm**

A straightforward way of solving the bilinear optimization problem (11) is to use a bilinear solver, whereas the bilinear solvers are not so developed as the linear ones [31]. A bypass solution is to relax the problem into linear subproblems and then use a two-way iterative search algorithm, which leads to the so-called V-s iteration algorithm [31,44]. Then the problem (11) turns out to be a three-linear problem with three steps, namely, $\gamma$-step, $\beta$-step, and $V$-step. When stopping criteria of the algorithm are met, a new LF is found with a proven level set

$$\Omega^* := \{x \in \mathbb{R}^n : V^* < 1\}, \tag{12}$$

where $V^*$ is the optimized LF with the largest possible level set as an estimate of the true ROA. Given an initial feasible LF $V_0$, a shape function $p_0$, and the expected number of iterations $N_I$, the V-s iteration algorithm, denoted as Algorithm 1, can be executed as follows

---

**Algorithm 1:** V-s **iteration**

**Input:** a proper LF $V_0(x)$; a shape function $p_0(x)$; the number of iterations $N_I$;
**Output:** $V^*$.
1:    $V = V_0$, $p = p_0$;
2:    **for** $i = 1:N_I$ **do**
3:        $\gamma$ **-step: hold** $V$ **fixed and solve for** $s_2$ **and** $\gamma^*$**:**

$$\gamma^* := \max_{s_2 \in SOS} \gamma \quad \text{s.t.} \quad -[(\partial V / \partial x)f + l_2] - (\gamma - V)s_2 \in \Sigma_n; \qquad (13)$$

4:        $\beta$ **-step: hold** $V$ **and** $\gamma^*$ **fixed and solve for** $s_1$ **and** $\beta^*$**:**

$$\beta^* := \max_{s_1 \in SOS} \beta \quad \text{s.t.} \quad (\gamma^* - V) - (\beta - p)s_1 \in \Sigma_n; \qquad (14)$$

5:        $V_{old} = V$, $\beta^*_{old} := \beta^*$;
6:        $V$ **-step: hold** $s_1, s_2, \beta^*, \gamma^*$ **fixed and solve for the new** $V$ **satisfying:**

$$\begin{aligned} -[(\partial V / \partial x)f + l_2] - (\gamma^* - V)s_2 &\in \Sigma_n \\ (\gamma^* - V) - (\beta^* - p)s_1 &\in \Sigma_n \\ V - l_1 &\in \Sigma_n \\ V(0) &= 0 \end{aligned} ; \qquad (15)$$

7:        **if** (15) **is feasible then**
            $V = V / \gamma^*$,
        **else**
            $V = V_{old} / \gamma^*$,
            **return**
        **end if**
8:        **if** $\left|(\beta^*_{old} - \beta^*) / \beta^*\right| < \varepsilon_{TOL}$ **then**
            **return**
        **end if**
9:    **end for**
10:   $V^* = V$, $\Omega^* := \{x \in \mathbb{R}^n : V^* < 1\}$.

---

**Remark 1:** To satisfy Eq.(14), the degree of SOS multiplier $s_1$ is chosen such that $\deg p + \deg s_1 \geq \deg V$. In Eqn.(13), $f(0) = 0$ indicates no constant term in $(\partial V / \partial x)f$ so that the multiplier $s_2$ associated with the term $(\gamma - V)$ in Eqn.(13) should not include a constant term either. Then the degree of $s_2$ is chosen larger than the maximum degree of $s_2$, $l_2$ and $(\partial V / \partial x) f$. For more details about the practical aspects of computation, readers are referred to [2,31].

**Remark 2:** The choice of an initial LF $V_0$ is flexible as long as it guarantees the initial search is feasible. Certainly, an optimized $V_0$ is favorable for better estimation, for instance, using the optimization procedure in [5] or simulation-guided procedure in [3] to find a better $V_0$. However, considering the difficulty of constructing LF, a more systematic way widely used in literature is

$$V_0 = x^T P x, \quad (16)$$

where $P$ is computed by the Lyapunov equation

$$A^T P + P A = -Q, \quad (17)$$

where $A = (\partial f / \partial x)\big|_{x=0}$ is considered Hurwitz. $P$ is a positive definite matrix. $Q > 0$ and $Q = I$ ($I$ is the unit matrix) is usually selected to provide the largest ball estimation [2].

**Remark 3:** In addition to expanding estimation interiorly, the shape function $p(x)$ can be taken as a size measurement of ROA. The convergence of $\beta$ and other elements can be taken as stopping criteria of the algorithm, such as

A. two consecutive $\beta$ is less than a programmed tolerance $\varepsilon_{TOL}$;

B. the numerical infeasibility alert in the optimization problem in Eqns. (13)-(15);

C. specified number of iterations $N_I$ is reached.

**Remark 4:** The V-s iteration algorithm converts the computation difficulty of a bilinear problem to an iteration procedure including two optimization problems ($\gamma$-step and $\beta$-step) and a feasibility problem ($V$-step) at each iteration that can be solved by linear SDP solvers, such as

SeDuMi. The bisection procedure has to be applied in $\gamma$-step and $\beta$-step because the optimization variable $\gamma$ ($\beta$) is coupled with the decision variables in $s_2$ ($s_1$).

**Remark 5:** The computation complexity grows dramatically with the scale of the problem, namely, the degree and dimension of the system $f(x)$ and the degree of LF being searched, which is the nature of SOS optimization. It roughly limits the approach to systems with fewer than eight to ten states with a cubic degree. Polynomial models of higher degrees can be handled with fewer states [47]. Usually, $\gamma$-step occupies the majority of the computation cost due to a higher overall degree. Software and additional information about V-s iteration can be found in [48]. The free distributed toolbox SOSOPTS [26] is used for examples in Sec. 5.

## 3 Preliminaries on shape function

It is mentioned in Sec. 2.1 that a user-defined shape function $p(x)$ is introduced to enlarge the inner approximation for the exact ROA and a proper choice of $p(x)$ can enhance the estimation [2]. According to the definition, $p(x)$ can be any positive definite and convex polynomial but a widely used choice is a quadratic function with its geometric center at the origin

$$p(x) = x^T N x, \tag{18}$$

where the positive definite matrix $N$ used for practical implementation. Reasons for choosing this form are three-fold. First, the introduction of $\beta$-step can incur additional computation costs. Therefore, the interior expansion step is implemented with a computationally cost-efficient quadratic form to avoid obscuring the estimation improvement by the incurred computation cost. Second, the system is assumed to have an asymptotically stable EP at the origin and thus the origin is the only identified point inside the true ROA without additional knowledge, which justifies the geometric center at the origin. Third, the shape matrix $N$ is vital for an accurate estimation because it can carry certain physical meanings, such as shape information of the ROA, dimensional scaling information as well as the importance of certain directions in the state space [44]. Its choice can be problem-dependent, but a general choice $N = I$ is made when prior knowledge is unavailable. Besides, the quadratic $V_0$ computed by Eqn. (16) can also be taken as a shape function.

In the traditional V-s iteration (**Algorithm 1**), the shape function is kept fixed throughout the iteration. A modification in [2] proposes to update the shape function iteratively by the quadratic part of a newly found LF. This modified algorithm (**Algorithm 2**) can be illustrated by an additional step after $V$-step in **Algorithm 1**.

| **Algorithm 2: Modified V-s iteration with an adaptive shape function** |
| --- |
| **1-7:**      Same with Step 1-7 in **Algorithm 1**; |
| **7a:**      Update the shape function $p(x)$ with the quadratic part of $V$ ; |
| **8-10:**      Same with Step 8-10 in **Algorithm 1**. |

Since $V$ has to be positive definite and $V(0) = 0$, it does not contain constant and linear terms. Therefore, the adaptive shape function will have and only have quadratic terms, which means that the geometric center of the shape function is still at the origin. As demonstrated in [2], the adaptive shape function aligns better with some simple-shaped ROA, for instance, ROA of the Van de Pol system, but struggles to do that for complex-shaped or unbounded ROA. In addition, it intends to align with only one direction (or converge in one direction) and possibly cause the estimation degradation in other directions. Meanwhile, it suffers from the convergence issue as **Algorithm 1** does.

For the sake of brevity, **Algorithm 1** and **Algorithm 2** will be designated by A1 and A2 respectively in the rest of the paper.

# 4 ROA estimation via RcomSSF

The before-mentioned advancements in estimation bring restricted enlargement with certain limitations, for example, the dramatic increase of computation burden, the early convergence, and only effective for certain systems or ROA. For the V-s iteration algorithm, with the geometric center of the shape function staying at the origin, the expansion is bounded in a domain around the origin thus limiting applications for non-symmetric and unbounded ROA. In this paper, the proposed RcomSSF improves ROA estimation by shifted shape function through shifting procedure and takes advantage of the existing advancements while avoiding their limitations. Besides, the advantage of leveraging the results already obtained will be demonstrated by the iteratively shifting in rounds of shifts.

## 4.1 Shifted Shape function

This section reveals the concept of the shifted shape function. Rather than fixed at the origin as the conventional shape function in Eqn. (18), the center of the shape function is shifted away from the origin to another point inside a verified level set. Since the center is a valid point inside the true ROA, the shifted shape function will not violate the definition of shape function and will still play the role of interior expansion and guide the estimation growth towards a new region. Besides, the shifted shape function is effective especially for a non-symmetric or unbounded ROA due to the offset of the shifting center. The shifts can be done to any shape function, but for simplicity in demonstration and implementation, an ellipsoidal shape function with a shifting center $x^*$ is considered in this paper

$$p(x) = (x - x^*)^T N (x - x^*). \tag{19}$$

It can be seen as an extension of Eqn. (18) with additional constant and linear terms brought by $x^*$. Though simple and can be absorbed in the general construction of SOS polynomials, it has never been explored, to the best knowledge of the authors. The benefits of using this form rather than a general SOS construction by monomials are that Eqn. (19) provides an efficient and systematic way of constructing shape functions through $x^*$ and $N$ that can carry explicit physical information and can be specified to improve ROA expansion capabilities. In addition, since the proposed shape function construction fits the SOS framework, the well-elaborated numerical optimization procedures facilitate our approach.

Shape functions with different shifting centers and shape matrices are illustrated by two-dimensional examples in Fig. 1. Comparison of $p_2(x)$ and $p_1(x)$ indicates the effect of different shape matrices, and comparison of $p_3(x)$ and $p_2(x)$ shows the shifting effect of $x^*$.

### 4.2 Rounds of shifts

The shifting procedure realized by rounds of shifts that leverage the results already obtained can be described below. A shape function with a new center $x^*$ will generally produce a new level set $\Omega^*$. Likewise, a newly obtained level set can continue to produce new sets by choosing new shifting centers inside. That is to say, a series of shifts will yield a series of different proven level sets by the iterative shifting. The process can be continued until convergence to the exact ROA. Afterward, these obtained level sets can be united into one set by R-composition (which will be explained later) [39,40] as the final verified inner approximation of the exact ROA, which is the proposed RcomSSF.

For a detailed description, a proven level set of the ROA has to be obtained first, for which V-s iteration can be used and the result is denoted as $\Omega_0^* := \{x \in \mathbb{R}^n : V_0^* < 1\}$. Algorithms A1, A2, or other advancements can be used to find a larger $\Omega_0^*$, which allows the RcomSSF to make full use of the existing advantageous algorithms. Then by shifting procedure, other proven level sets will stem from $\Omega_0^*$.

As shown in Fig. 2, a set of shifting centers $x_i^*$ ($i = 1, 2, ...,$ $i$ indexes the shifts in the first round) are chosen from $\Omega_0^*$ to form the first round of shifts. Then the V-s algorithm is applied to obtain the corresponding optimized LFs $V_i^*$ and the relevant proven level sets $\Omega_i^* := \{x \in \mathbb{R}^n : V_i^* < 1\}$. Next, for the second round of shifts, centers $x_{ij}^*$ ($j = 1, 2, ...$ indexes the shifts in the second round) are chosen inside the respective proven level sets $\Omega_i^*$, and then produce the second generation of optimized LFs $V_{ij}^*$ and the corresponding proven level sets $\Omega_{ij}^* := \{x \in \mathbb{R}^n : V_{ij}^* < 1\}$. Similarly, the third round of shifts is done to obtain the third generation of proven level sets $\Omega_{ijk}^* := \{x \in \mathbb{R}^n : V_{ijk}^* < 1\}$ ($k = 1, 2, ...$ indexes the shifts in the third round). Ideally, the process is in a tree structure as seen in Fig. 2. To differentiate different shifts, the indexes $i, j, k$ denote the independent shifts in the first, second and third round, respectively. The subscripts, except the last one, refer to the parent LF and the parent level set, for instance, $\Omega_{123}^* := \{x \in \mathbb{R}^n : V_{123}^* < 1\}$ is obtained by the center $x_{123}^*$ chosen from the parent level set $\Omega_{12}^* := \{x \in \mathbb{R}^n : V_{12}^* < 1\}$. For the sake of brevity, the approach is demonstrated only for three rounds of shifts, however, without loss of generality, it can be extended to a higher number of rounds.

It needs to be mentioned that the V-s iteration algorithm dealing with the shifted shape function must have a fixed center at the shifting center to maintain the expansion around the specific center. Otherwise, the shifting effect will be neutralized, for example, if A2 is used, the adaptive shape function composed by the quadratic terms of the new LF will always bring the center back to the origin. Examples in Section 5 will show how this works.

### 4.3 Selection of shifting centers

Here some considerations for selecting the shifting center $x^*$ are provided.

Firstly, $x^*$ must be contained inside a proven level set $\Omega$ to guarantee the set containment constraint in Eqn. (7) is satisfied thus ensuring the following iterations are feasible.

Secondly, $x^*$ should be close to the boundary of $\Omega$. For the sake of clarity, a two-dimensional example is used for demonstration in Fig. 3. Several locations (denoted by $C_i$, $i = 1, 2, ...$) of $x^*$ are shown inside $\Omega$. The location $C_1$ is expected to get more proved region outside $\Omega$ than $C_2$ because $C_1$ is closer to the boundary. One way is to choose the location of $x^*$ manually after getting a visualized proven level set; while for a general case, we propose a programmable procedure to assist the selection. As shown in Fig. 3, given a phase angle $\alpha$, calculate the distance $\rho_\alpha$ from the origin to the boundary, and determine the intersection point $P_\alpha$. Then the shifting center can be chosen at $x^* = \sigma \cdot x_{P_\alpha}$, where $x_{P_\alpha}$ is the coordinate of $P_\alpha$ and $\sigma \in (0,1)$. Simulation shows that for a smaller value of $\sigma$, a high number of iterations is required for the estimation to grow beyond the parent level set; otherwise, for a large value of $\sigma$, the center is too close to the boundary, and the proximity means a small $\beta$ in Eqn. (6) at the beginning and cannot expand the estimation properly. Therefore, a tradeoff has to be done to opt for a better location and a recommended value of $\sigma$ around 0.8 is given. The phase angle $\alpha$ specifies the growing direction of the estimation so that it can be decided according to the area that requires verification.

Thirdly, a center $x^*$ near the convex boundary usually yields better estimation. On the contrary, a center near the concave boundary sees little estimation expansion because usually, it has already approached the concave boundary of the exact ROA. For example, $C_i(i = 1, 4, 5, 6)$ are better than $C_3$.

Fourthly, there is no limit on the number of selected centers and, in general, the higher number the more accurate estimation. However, one has to note that increasing the number of selected centers also increases the computation cost. Therefore, a tradeoff is suggested between accuracy and computation cost so as to choose a proper number of centers.

These considerations for shifting center selection will help in understanding the method and Sec. 5 will show how these considerations work for specific examples.

### 4.4 R-composition

R-composition, a systematic way of composing LFs through the use of R-functions, has been proposed to obtain richer and more flexible LFs and used for ROA estimation in [39,40]. Here, it is employed in the final step of RcomSSF to provide a more compact result by uniting these independent level sets obtained from the shifting procedure into a single one. The compact form of the result is also beneficial for further applications of the result. It thus makes the proposed method more complete.

R-functions are important in R-composition because they represent the natural extension of Boolean operators (e.g. AND, OR, NOT) to real-valued functions and thus provide the basic tools to compute an analytic expression of intersections union, and complement in a geometric setting because. The full account of R-functions goes beyond this work and can be found in [39] and references therein. Constructed by the LF $V^*(x)$, the function

$$R(x) = 1 - V^*(x) \tag{20}$$

is an R-function. When $R(x) > 0$, the set $\hat{R} = \{x \in \mathbb{R}^n : R(x) > 0\}$ is exactly the proven level set in Eqn. (12) for ROA estimation. Possible choices, R-negation, R-intersection, and R-union, are outlined in Table 1.

The parameter $\tau$ chosen within (0, 2] keeps some implementation freedom and in this work $\tau = 2$ is chosen. Geometrically, if $R_1(x)$ and $R_2(x)$ are positive inside a geometrical region and negative outside, then the R-intersection and R-union represent the intersection and union between the two sets $\hat{R}_1 := \{x \in \mathbb{R}^n : R_1(x) > 0\}$ and $\hat{R}_2 := \{x \in \mathbb{R}^n : R_2(x) > 0\}$, which is

$$\hat{R}_1 \cap \hat{R}_2 = \{x \in \mathbb{R}^n : R_\cap(R_1, R_2) > 0\}, \quad \hat{R}_1 \cup \hat{R}_2 = \{x \in \mathbb{R}^n : R_\cup(R_1, R_2) > 0\}. \tag{21}$$

For instance, if $R_1(x) = 1 - x^T N_1 x$ and $R_2(x) = 1 - x^T N_2 x$, where $N_1 = \text{diag}(1,9)$ and $N_2 = \text{diag}(9,1)$, then the intersection and the union between sets $\hat{R}_1$ and $\hat{R}_2$ are shown in 错误!未找到引用源。.

Therefore, R-union is adopted in RcomSSF to produce the union of estimation. For the example of rounds of shifts presented in Fig. 2 with the obtained LFs $V_i^*, V_{ij}^*, V_{ijk}^*, ...(i, j, k = 1, 2, 3...)$, the final estimation of ROA can be computed iteratively by

$$\Omega_e = \{x \in \mathbb{R}^n : R_e > 0\}, \quad R_e = R_\cup(R_\cup(R_\cup(1-V_0^*, 1-V_i^*), 1-V_{ij}^*), 1-V_{ijk}^*, ...) \tag{22}$$

and we have

$$\Omega_e = \Omega_0^* \cup \Omega_i^* \cup \Omega_{ij}^* \cup \Omega_{ijk}^* \cdots. \tag{23}$$

It has been discussed in [39,40] that $-R_e(x) + R_e(0)$ is a Lyapunov function when $\tau = 2$ or more precisely, a Lyapunov-like function since the classical Lyapunov function condition that requires continuous differentiability is relaxed. Sufficient conditions for $-R_e(x) + R_e(0)$ being a LF can be found in [39,40] and a detailed discussion is omitted here since R-composition in RcomSSF takes on the task of better communicating the result by a compact form and the union of LFs is not used to initialize the next optimization process though the possibility of doing this is held.

## 4.5 RcomSSF

Based on the details above, the algorithm can be summarized as follows:

**Algorithm 3: RcomSSF**

**Input:** a proper LF $V_0(x)$; shape function $p_0(x)$; the number of iteration $N_I$

**Output:** $\Omega^*$

1: Run V-s iteration to obtain the optimized LF $V_0^*$ and the relevant parent level set $\Omega_0^* := \{x \in \mathbb{R}^n : V_0^* < 1\}$;

2: The first round of shifts:

   **Initialization:**
   Give phase angle $\alpha$ and coefficient $\sigma$ to calculate $\rho_\alpha$ and locate shifting centers $x_i^*$ ($i = 1, 2, ...$) in $\Omega_0^*$; specify the shape matrix $N$; $V_0 = V_0^*$, $p_0 = (x - x_i^*)^T N(x - x_i^*)$ and iteration number $N_I$;

   **Run V-s iteration:** obtain the optimized LF $V_i^*$ and level set $\Omega_i^*$;

   **Further shift check:**
   Calculate $\rho_\alpha$ again as $\rho_{\alpha New}$
   **if** $|(\rho_{\alpha New} - \rho_\alpha) / \rho_\alpha| > \varepsilon_{TOL\rho}$ **then**
       continue the next round;
   **else**
       go to Step 6;
   **end if**

3: The second round of shifts:

   **Initialization:**
   Give phase angle $\alpha$ and coefficient $\sigma$ to calculate $\rho_\alpha$ and locate shifting centers $x_{ij}^*$ ($j = 1, 2, ...$) in $\Omega_i^* := \{x \in \mathbb{R}^n : V_i^* < 1\}$; specify the shape matrix $N$; $V_0 = V_i^*$, $p_0 = (x - x_{ij}^*)^T N(x - x_{ij}^*)$ and iteration number $N_I$;

   **Run V-s iteration:** obtain the optimized LF $V_{ij}^*$ and level set $\Omega_{ij}^*$

   **Further shift check.**

4: The third round of shifts:

   **Initialization:**
   Give phase angle $\alpha$ and coefficient $\sigma$ to calculate $\rho_\alpha$ and locate shifting centers $x_{ijk}^*$ ($k = 1, 2, ...$) in $\Omega_{ij}^* := \{x \in \mathbb{R}^n : V_{ij}^* < 1\}$; specify the shape matrix $N$; $V_0 = V_{ij}^*$, $p_0 = (x - x_{ijk}^*)^T N(x - x_{ijk}^*)$ and iteration number $N_I$;

   **Run V-s iteration (A1):** obtain the optimized LF $V_{ijk}^*$ and level set $\Omega_{ijk}^*$;

   **Further shift check.**

| 5: | The fourth round of shifts … |
|---|---|
| 6: | R-composition of $\Omega_0^*$, $\Omega_i^*$, $\Omega_{ij}^*$, and $\Omega_{ijk}^*$ …into one single level set $\Omega_e$ by Eqn. (22). |

The proximity of the center to the boundary $\left|(\rho_{\alpha New} - \rho_\alpha)/\rho_\alpha\right|$ is taken as a metric for further shift check. The tolerance $\varepsilon_{TOL\rho}$ can be customized, such as 10%, which means if the expansion in phase $\alpha$ direction is less than 10% then the shifting stops in this direction; otherwise, a new center will be decided and the next round of shifts will be carried out. And shifts in the same round can be carried out in parallel to shorten the whole verification period.

In RcomSSF, the increase of shift rounds with an increased number of selected shift centers leads to a linear increase in the number of iterations. It brings the advantage of RcomSSF in computation cost over some other modifications, for example, improvement by using an increased degree of LF. Moreover, the shifting procedure yields more significant growth in estimation. These considerations articulate the computational efficiency of the proposed method. Besides, when applying to practical problems, RcomSSF is feasible as long as the V-s iteration algorithm is feasible and RcomSSF yields better estimation without other more sophisticated approaches.

Generally, prior knowledge of the true ROA is favorable in the estimating process, however, is not a must for RcomSSF. RcomSSF is effective even without searching for special shape functions or locations of the shifting center. Certainly, existing optimization techniques for the initial LF and shape matrix from literature (e.g. [3,5,40,44]) can still be absorbed in RcomSSF and work for improvement.

In the next section, RcomSSF will be compared on several benchmark examples with A1, A2, and some other methods from literature to show its effectiveness and advantages.

# 5 Examples

## 5.1 Systems with bounded ROAs

**Example 1.** Here we consider a Van de Pol system taken from [2]

$$\begin{cases} \dot{x}_1 = -x_2 \\ \dot{x}_2 = x_1 + 5x_2(x_1^2 - 1) \end{cases}. \tag{24}$$

It has a stable equilibrium point (EP) at the origin and an unstable limit cycle. The problem of finding the ROA of Van de Pol systems has been studied extensively [2,3,19,38,49] and it is thus taken as a benchmark example for testing the new ROA estimation method. Its exact ROA is the region enclosed by the limit cycle and can be plotted by the reverse trajectory method as shown in Fig. 5 (a). To start the estimation, an initial LF $V_0(x)$ is computed by Eqn. (16) 错误!未找到引用源。 using $Q = I$ and then matrices $A$ and $P$ are

$$A = \begin{bmatrix} 0 & -1 \\ 1 & -5 \end{bmatrix}, \quad P = \begin{bmatrix} 2.7 & -0.5 \\ -0.5 & 0.2 \end{bmatrix}.$$

Then the computation is initialized by

$$V_0(x) = 2.7x_1^2 - x_1 x_2 + 0.2 x_2^2, \quad p_0(x) = V_0(x). \tag{25}$$

Estimation results presented in the recent work [2] using algorithms A1 and A2 under the same initial conditions are also shown for comparison with the proposed RcomSSF. It can be seen in Fig. 5 (a) that when a six-degree LF is searched with the iteration number $N_I = 30$ and 60, A2 with an adaptive shape function yields a larger estimation than A1 with a fixed shape function but still fails to occupy the whole exact ROA. In fact, there is no improvement in estimation for both A1 and A2 after 20 iterations because the algorithm converges. This is elaborated in further detail for A2 in Fig. 5 (b). Variation of the adaptive shape function is plotted for the 1st, 30th, and 60th iteration. It shows that the ellipse rotates anticlockwise to align better with the shape of the true

ROA. Meanwhile, stretching along the major axis and shrinking along the minor axis can also be observed, which subsequently leads to a slight increase and decrease of estimation in the corresponding direction. However, the meaningful rotating and stretching almost stop after 30 iterations, and correspondingly the estimation sees convergence. In summary, A1 does not perform as well as A2 for the same amount of computation (same number of iterations); both A1 and A2 have the limitation of early convergence and thus a further increase of iterations does not help in expanding the estimation; the shrinking of shape function in a certain direction gives degradation of the estimation.

Now the proposed RcomSSF is utilized to overcome the above limitations. The V-s iteration is performed first to yield a valid initial level set. The final $V^*$ from A1 at $N_I = 30$ is selected as the initial LF $V_0^*$. A shape matrix $N = I$ is selected for a general circle $p(x)$. To initialize the first round of shifts, two centers $x_1^* = [1,1]$ and $x_2^* = [-1,-1]$ are picked directionally inside the level set $\Omega_0^* := \{x \in \mathbb{R}^n : V_0^* < 1\}$ given the gap between the exact ROA. Subsequently, the shape function is constructed by Eqn. (19) with $x_1^* = [1,1]$, which is

$$p_1(x) = x_1^2 + x_2^2 - 2x_1 - 2x_2 + 2.$$

By RcomSSF, a proven level set $\Omega_1^* := \{x \in \mathbb{R}^n : V_1^* < 1\}$ (green dotted line in Fig. 5 (c)) is thus obtained, which fills the gap on the right corner nicely. For the center $x_2^* = [-1,-1]$,

$$p_2(x) = x_1^2 + x_2^2 + 2x_1 + 2x_2 + 2.$$

Likewise, a proven level set $\Omega_2^* := \{x \in \mathbb{R}^n : V_2^* < 1\}$ (green dashed line in Fig. 5 (c)) is obtained, which fills the gap on the left corner nicely. Together, the union of the three level sets $\Omega_0^*$, $\Omega_1^*$, and $\Omega_2^*$ is computed by R-composition by Eqn. (22)

$$\Omega_e := \{x \in \mathbb{R}^n \mid R_e > 0\}, \quad R_e = R_U(R_U(1-V_0^*, 1-V_1^*), 1-V_2^*) \tag{26}$$

and

$$\Omega_e = \Omega_0^* \cup \Omega_1^* \cup \Omega_2^*. \tag{27}$$

As can be seen from Table 1, $R_e$ is a complex function. However, the polynomial approximation can be found and gives the following representation

$$\begin{aligned}
R_e(x) = &-0.097268 x_1^6 - 0.047707 x_1^5 x_2 - 0.007790 x_1^4 x_2^2 + 0.003634 x_1^3 x_2^3 \\
&- 0.001048 x_1^2 x_2^4 + 0.000208 x_1 x_2^5 - 3.601840 \mathrm{e}^{-5} x_2^6 + 0.000273 x_1^5 \\
&- 0.000271 x_1^4 x_2 + 0.000396 x_1^3 x_2^2 - 9.588542 \mathrm{e}^{-5} x_1^2 x_2^3 + 1.126447 \mathrm{e}^{-5} x_1 x_2^4 \\
&- 2.250610 \mathrm{e}^{-7} x_2^5 + 0.542801 x_1^4 + 0.073756 x_1^3 x_2 - 0.082728 x_1^2 x_2^2 \\
&+ 0.021571 x_1 x_2^3 - 0.003026 x_2^4 - 0.000389 x_1^3 - 0.000149 x_1^2 x_2 \\
&- 0.000249 x_1 x_2^2 + 8.912165 \mathrm{e}^{-5} x_2^3 - 1.182416 x_1^2 + 0.417432 x_1 x_2 \\
&- 0.076201 x_2^2 + 1
\end{aligned}.$$

Another comparison from [45] is discussed here. A numerical method is proposed in [45] with the conjectured convergence to the true ROA arbitrarily well. However, the result obtained with a six-degree LF still presents an obvious gap with the true ROA for the Van de Pol example, similar to the result of A2 shown in Fig. 5 (a). Furthermore, using a higher-degree LF makes this method almost inapplicable for real-world problems due to the dramatic growth of computational burden. Besides, the method [45] is limited only to bounded ROA but the proposed RcomSSF also applies to unbounded ROA as illustrated in the following examples.

## 5.2 Systems with unbounded ROAs

**Example 2.** Consider the following system [2,50]

$$\begin{cases} \dot{x}_1 = -4x_1^3 + 6x_1^2 - 2x_1, \\ \dot{x}_2 = -2x_2. \end{cases} \tag{28}$$

Analysis of the linearized system, combined with the vector field plot in Fig. 6 (a), shows that the system has three EP, namely, two stable node sinks (0,0) and (1,0) and the saddle point (0.5, 0). A line $x_1 = 0.5$ divides the plane. The trajectories originate from the left-hand/ right-hand side of the line sink down to node (0, 0) /(1,0), thus two ROA corresponding to the two stable node sinks. Each of them can be treated equivalently, however, only EP (0,0), whose ROA is unbounded in $x_1 < 0.5$ plane, is considered in this study to evaluate the technical performance of RcomSSF.

For ROA estimation, the initial LF is computed by Eqn. (16) using $Q = I$, and matrices $A$ and $P$ are obtained

$$A = \begin{bmatrix} -2 & 0 \\ 0 & -2 \end{bmatrix}, \quad P = \begin{bmatrix} 0.25 & 0 \\ 0 & 0.25 \end{bmatrix}.$$

Then the computation is initialized by

$$V_0(x) = 0.25 x_1^2 + 0.25 x_2^2. \tag{29}$$

We will first look at the results of algorithms A1 and A2 [2,50]. For them, the initial shape function is set as

$$p_0(x) = 0.8 \cdot V_0(x) \tag{30}$$

and a quartic LF is searched. As shown in Fig. 6 (a), A1 converges after about 30 iterations, but A2 gets a larger estimation after the same number of iterations by updating the shape function. The adaptive shape function takes the form of

$$p(x) = 2.58 x_1^2 + 3.47 \mathrm{e}^{-2} x_2^2 - 3.28 \mathrm{e}^{-8} x_1 x_2 \tag{31}$$

at the 30$^{\text{th}}$ iteration and covers a much larger area as shown in Fig. 6 (a). But simulation result in [2] shows that the left boundary of $x_1$ is still limited by $-1.0$ even when iterations are increased to 150. This limitation, however, can be lifted by the proposed RcomSSF that allows an extra extension in $-x_1$ direction using a shifted shape function. The level set obtained by A2 $\Omega_0^* := \{x \in \mathbb{R}^n : V_0^* < 1\}$ for 30 iterations is used to produce shifting centers and a center $x_1^* = [-0.8, 0]$ is chosen to expand towards the left. With the shape matrix $N = I$, the shifted shape function constructed by Eqn. (19) is

$$p_1(x) = x_1^2 + x_2^2 + 1.60x_1 + 0.64. \tag{32}$$

Using RcomSSF, a new proven level set $\Omega_1^* := \{x \in \mathbb{R}^n : V_1^* < 1\}$ (green dotted line in Fig. 6 (b)) is obtained. Then the left boundary is extended from $-1.0$ to $-1.8$ but the estimation along $x_2$ axis shrinks. Inspired by the adaptive shape function in A2, a shape matrix $N_1 = \mathrm{diag}(1, 1/16)$ is customized in order to put more weight on $x_2$ axis. In this case,

$$p_{1_{N_1}}(x) = x_1^2 + 0.0625 x_2^2 + 1.60 x_1 + 0.64. \tag{33}$$

As expected, $N_1$ produces a larger level set after 30 iterations (green dash-dot line) and after 60 iterations (green dashed line, denoted as $\Omega_{1_{N_1}}^* := \{x \in \mathbb{R}^n : V_{1_{N_1}}^* < 1\}$). The estimation expands not only along $-x_1$ axis but also $x_2$ axis. The union of those obtained level sets is computed by Eqn. (22) as

$$\Omega_e := \{x \in \mathbb{R}^n \mid R_e > 0\}, \quad R_e = R_\cup(1 - V_0^*, 1 - V_{1_{N_1}}^*) \tag{34}$$

and

$$\Omega_e = \Omega_0^* \cup \Omega_{1_{N_1}}^*. \tag{35}$$

The resulting level set is represented by the red dotted line in Fig. 6 (b).

This example shows the ability of the proposed RcomSSF method in extending the estimation in a specific direction.

**Example 3.** Consider the following system [2]

$$\begin{cases} \dot{x}_1 = -50x_1 - 16x_2 + 13.8x_1x_2 \\ \dot{x}_2 = 13x_1 - 9x_2 + 5.5x_1x_2 \end{cases}. \tag{36}$$

Analysis of the linearized system gives a stable node (0, 0) and a saddle point (1.45, 18.17). The vector field plot shows the boundary of the ROA for EP (0, 0) in Fig. 7 (a). At the initialization stage of the algorithm, $V_0(x)$ computed by Eqn. (16) with $Q = I$ is

$$V_0(x) = 0.011694x_1^2 + 0.013034x_1x_2 + 0.043969x_2^2. \tag{37}$$

The shape function is set $p_0(x) = V_0(x)$ and a four-degree LF is searched iteratively.

In Fig. 7 (a), the results of algorithms A1 (blue lines) and A2 (black lines) are shown. The fact that there is no significant improvement after 30 iterations indicates convergence is achieved. It again reinforces the idea that the increase of iteration number cannot be used as a universal remedy. Besides, estimation obtained by A1 approaches the exact boundary only in a certain region, and the adaptive shape function in A2 guides the estimation to grow in a skew direction and does not fit the boundary well. However, the half-unbounded ROA means that there is still plenty of room for the estimation to grow. The proposed RcomSSF realizes the extension by the shifting procedure. The union of two rounds of shifts is plotted in Fig. 7 (a) as well for comparison purpose. Effects of selecting a set of shifting centers and specifying the shape matrix $N$ are demonstrated in Fig. 7 (b) and (c).

To start the first round of shifts, the proven level set at the 30$^{th}$ iteration by A1 is chosen as the parent level set $\Omega_0^* := \{x \in \mathbb{R}^n : V_0^* < 1\}$ to produce shifting centers. According to the considerations given in Sec. 4.3, two centers $x_1^* = [0, -4]$ and $x_2^* = [-7.5, 0]$ are chosen to pull the estimation down and left. First, for $x_1^* = [0, -4]$, the shape matrices $N_1 = \mathrm{diag}(1/4, 1)$ is attempted for a larger expansion and the corresponding shape function is

$$p_1(x) = 0.25 x_1^2 + x_2^2 + 8x_2 + 16. \tag{38}$$

The obtained proven level sets are shown by green lines (green solid line $\Omega_1^* := \{x \in \mathbb{R}^n : V_1^* < 1\}$ and green dashed line $\Omega_2^* := \{x \in \mathbb{R}^n : V_2^* < 1\}$) in Fig. 7 (b), which present significant extension. Then the second round of shifts is carried out and the proven level sets are shown by magenta lines in Fig. 7 (b). The center $x_{11}^* = [0, -11]$ is selected from the previously obtained $\Omega_1^*$ to move the extension further downwards. With a shape matrix $N_1 = \mathrm{diag}(1/4, 1)$, the shifted shape function is

$$p_{11}(x) = 0.25 x_1^2 + x_2^2 + 22 x_2 + 121 \tag{39}$$

and the resulting level set $\Omega_{11}^* := \{x \in \mathbb{R}^n : V_{11}^* < 1\}$ is represented by the magenta solid line in Fig. 7 (b). Then in $\Omega_2^*$, centers $x_{21}^* = [-18, 2]$ and $x_{22}^* = [-3, 8]$ are chosen to pull the estimation towards the left and top, respectively. With $N = I$, the shifted shape functions are the following

$$p_{21}(x) = x_1^2 + x_2^2 + 36 x_1 - 4 x_2 + 328, \quad p_{22}(x) = x_1^2 + x_2^2 + 6 x_1 - 16 x_2 + 73. \tag{40}$$

The resulting level sets are represented by a magenta dashed line ($\Omega_{21}^* := \{x \in \mathbb{R}^n : V_{21}^* < 1\}$) and magenta dotted line ($\Omega_{22}^* := \{x \in \mathbb{R}^n : V_{22}^* < 1\}$) in Fig. 7 (b). Given the limited space, further shifts will be omitted here. Finally, the union of the level set is computed by Eqn. (22) as

$$\begin{aligned} &\Omega_e := \{x \in \mathbb{R}^n \mid R_e > 0\}, \\ &R_e = R_\cup(R_\cup(R_\cup(R_\cup(R_\cup(1 - V_0^*, 1 - V_1^*), 1 - V_2^*), 1 - V_{11}^*), 1 - V_{21}^*), 1 - V_{22}^*) \end{aligned} \tag{41}$$

and

$$\Omega_e = \Omega_0^* \cup \Omega_1^* \cup \Omega_2^* \cup \Omega_{11}^* \cup \Omega_{21}^* \cup \Omega_{22}^*. \tag{42}$$

The union $\Omega_e$ can be visualized in Fig. 7 (a).

From the comparison of A1, A2, and RcomSSF methods shown in Fig.7(a), it could be concluded that early convergence occurs for A1 and A2 thus an increased number of iterations (30 to 60) cannot bring further expansion, which is improved by the proposed RcomSSF that manifests significant estimation improvement after only several shifts.

One may note that $N_1 = \text{diag}(1/4,1)$ is chosen rather than $N = I$ for the two shifts downwards. It is elaborated in Fig. 7 (c) where estimation obtained by $N_1$ expands more along $x_2$ axis since more weight is put on $x_2$ in $N_1$. This gives hint on the effect of the shape matrix and its selection.

**Example 4.** Here we consider the Hahn example

$$\begin{cases} \dot{x}_1 = -x_1 + 2x_1^2 x_2 \\ \dot{x}_2 = -x_2 \end{cases}. \tag{43}$$

This system has an asymptotically stable EP at the origin and its exact ROA is known as $x_1 x_2 < 1$. The determination of its ROA has been studied extensively [19,31,51]. A predetermined shape function

$$p(x) = x^T \begin{bmatrix} 14.47 & 18.55 \\ 18.55 & 26.53 \end{bmatrix} x \tag{44}$$

was used in [31]. Also in [31], the ROA estimation can be increased with higher-degree LFs and a composed LF by pointwise maximum or minimum of polynomials. Even though, there still exist uncovered regions near the origin. To further cover the remaining regions as close to the stability region $x_1 x_2 < 1$ as possible, multiple shape functions obtained by rotating the major axis of an ellipse every three degrees are attempted in [31] as well and a series of level sets are obtained. The envelope of these sets is shown in Fig. 8 (a) for comparison purpose. Though improved compared with the single LF method, the estimation is still boxed inside $|x| < 6$. Besides, a recent result obtained by an invariant set method [51] appears smaller even than that in [31]. Overall, these methods provide insufficient performance due to the unbounded nature of the problem.

For comparison, algorithms A1, A2, and the proposed RcomSSF are carried out with an initial $V_0(x)$ computed by Eqn. (16) with $Q = I$

$$V_0(x) = 0.5x_1^2 + 0.5x_2^2. \tag{45}$$

The initial shape function is set as Eqn. (44) for comparison with [31] and a six-degree LF is searched. As shown in Fig. 8 (a), the adaptive shape function in A2 underperforms the fixed shape function in A1, which indicates that the modification in A2 is not necessarily effective. Besides, both A1 and A2 face the problem of early convergence after certain iterations. By contrast, the proposed RcomSSF outperforms the methods mentioned above by extending the estimation significantly and expanding closely along the boundary. It is achieved by two rounds of shifts as shown in Fig. 8 (b). To start RcomSSF, the level set $\Omega_0^* := \{x \in \mathbb{R}^n : V_0^* < 1\}$ obtained by A1 is taken as the parent level set. Shifting centers of the first round are selected from $\Omega_0^*$ according to Sec. 4.3. Considering that the exact boundary $x_1 x_2 < 1$ and $\Omega_0^*$ are symmetric about the origin, two centers $x_1^* = [-4, 3]$ and $x_2^* = [4, -3]$ are chosen to produce the relevant shape functions by Eqn. (19) with $N = I$,

$$p_1(x) = x_1^2 + x_2^2 - 8x_1 + 6x_2 + 25 \text{ and } p_2(x) = x_1^2 + x_2^2 + 8x_1 - 6x_2 + 25. \tag{46}$$

Two proven level sets are thus obtained and shown by green lines in Fig. 8 (b) (green solid line $\Omega_1^* := \{x \in \mathbb{R}^n : V_1^* < 1\}$ and green dashed line $\Omega_2^* := \{x \in \mathbb{R}^n : V_2^* < 1\}$. Next, for the second round of shifts, centers $x_{11}^* = [-5, 5]$ and $x_{12}^* = [-6, 2]$ are selected from $\Omega_1^*$ and symmetrical centers $x_{21}^* = [5, -5]$ and $x_{22}^* = [6, -2]$ are selected from $\Omega_2^*$. The resulting level sets are shown by cyan lines accordingly. The second round approaches the upper and lower boundary closely. Moreover, the proposed routine can still be performed for the next rounds of shifts for further expansion. By R-composition in Eqn. (22), these independent level sets are united into a single set

$$\Omega_e := \{x \in \mathbb{R}^n \mid R_e > 0\},$$
$$R_e = R_\cup(R_\cup(R_\cup(R_\cup(R_\cup(1-V_0^*, 1-V_1^*), 1-V_2^*), 1-V_{11}^*), 1-V_{12}^*), 1-V_{21}^*), 1-V_{22}^*)$$
(47)

and

$$\Omega_e = \Omega_0^* \cup \Omega_1^* \cup \Omega_2^* \cup \Omega_{11}^* \cup \Omega_{12}^* \cup \Omega_{21}^* \cup \Omega_{22}^*$$
(48)

as shown in Fig. 8 (a).

**Example 5.** This example is a three-degree Taylor expansion of the system given in [52]

$$\begin{cases} \dot{x}_1 = x_2 + x_3^2 \\ \dot{x}_2 = x_3 - x_1^2 - x_1(x_1 - 1/6 x_1^3) \\ \dot{x}_3 = -x_1 - 2x_2 - x_3 + x_2^3 + 1/10(2/3 x_3^3 + 2/5 x_3^5) \end{cases}.$$
(49)

It has the origin as an asymptotically stable EP (0,0) and three other unstable EPs. ROA estimation results of EP (0,0) by A1, A2, and the proposed RcomSSF method are presented in Fig. 9 (a) and (b); cross-sections for $x_2 = 0$ are given in Fig. 9 (c); the three-dimensional visualization is given in Fig. 9 (d). It can be seen that the RcomSSF gives a larger provable level set by shifting to centers $x_1^* = [0.8, 0, 0]$, $x_2^* = [-0.8, 0, 0.6]$, $x_3^* = [0.2, 0, -0.8]$ and $x_4^* = [0, 0, 1.2]$ that are chosen in the level set obtained by A1. Denote $\Omega_0^*$ the set obtained by A1, $\Omega_i^*(i=1,2,3,4)$ the sets obtained in the shifting procedure, their union by R-composition using Eqn. (22) is

$$\Omega_e := \{x \in \mathbb{R}^n \mid R_e > 0\},$$
$$R_e = R_\cup(R_\cup(R_\cup(R_\cup(1-V_0^*, 1-V_1^*), 1-V_2^*), 1-V_3^*), 1-V_4^*)$$
(50)

and

$$\Omega_e = \Omega_0^* \cup \Omega_1^* \cup \Omega_2^* \cup \Omega_3^* \cup \Omega_4^*.$$
(51)

It certificates the effectiveness of RcomSSF on a higher dimensional system.

# 6 Conclusion

Knowledge of ROA is crucial for the analysis of nonlinear systems and control design. For a nonlinear polynomial system, this paper proposes a general method (RcomSSF) for ROA estimation improvement via an algorithm using shifted shape functions with centers shifted iteratively close to the boundary of the newly obtained proven subset of ROA. The algorithm is generally effective, even for non-symmetric or unbounded ROA, for which the existing methods present limitations. A composition method for Lyapunov functions, namely R-composition, is used in the proposed RcomSSF to unite the resulting independent level sets into a single level set or a single function, which brings a compact and richer-shaped expression.

Simulation results have been obtained for five benchmark examples from literature, including two- and three-dimensional systems with bounded or unbounded, symmetric, or non-symmetric ROA. Compared with the existing ROA estimation algorithms, the exceptional performance of RcomSSF is highlighted. It has the advantages, including improving estimation even when other algorithms encounter early convergence or numerical infeasibility; guiding the estimation into a required domain by specifying shifting centers or shape matrices; improving estimation only by a linear increase in computation burden instead of the nonlinear growth in the case of using high-degree LFs; maximizing the capability of existing advancements for ROA estimation; leveraging the results already obtained; being compatible with the existing optimization techniques for initial Lyapunov function and shape matrix; requiring no prior knowledge of the exact ROA or empirical assistance. It demonstrates an avenue for effective estimation of ROA that allows such implementation for real-world problems.

# Acknowledgments

The authors would like to acknowledge the time and effort devoted by reviewers to improving the quality of this work.

# References

[1] Khalil HK. Nonlinear systems. 3rd ed. Prentice Hall; 2002.
[2] Khodadadi L, Samadi B, Khaloozadeh H. Estimation of region of attraction for polynomial nonlinear systems: a numerical method. ISA Trans 2014;53:25–32. https://doi.org/10.1016/j.isatra.2013.08.005.
[3] Topcu U, Packard A, Seiler P. Local stability analysis using simulations and sum-of-squares programming. Automatica 2008;44:2669–75. https://doi.org/10.1016/j.automatica.2008.03.010.
[4] Yuan G, Li Y. Estimation of the regions of attraction for autonomous nonlinear systems. Trans Inst Meas Control 2019;41:97–106. https://doi.org/10.1177/0142331217752799.
[5] Sidorov E, Zacksenhouse M. Lyapunov based estimation of the basin of attraction of Poincare maps with applications to limit cycle walking. Nonlinear Anal Hybrid Syst 2019;33:179–94. https://doi.org/10.1016/j.nahs.2019.03.002.
[6] Cunis T, Condomines J-P, Burlion L. Local stability analysis for large polynomial spline systems. Automatica 2020;113:1–5. https://doi.org/10.1016/j.automatica.2019.108773.
[7] Giesl P, Osborne C, Hafstein S. Automatic determination of connected sublevel sets of CPA Lyapunov functions. SIAM J Appl Dyn Syst 2020;19:1029–56. https://doi.org/10.1137/19M1262528.
[8] Cunis T, Condomines J-P, Burlion L. Sum-of-squares flight control synthesis for deep-stall recovery. J Guid Control Dyn 2020;43:1498–511. https://doi.org/10.2514/1.G004753.
[9] Riah R, Fiacchini M, Alamir M. Iterative method for estimating the robust domains of attraction of non-linear systems: application to cancer chemotherapy model with parametric uncertainties. Eur J Control 2019;47:64–73. https://doi.org/10.1016/j.ejcon.2018.12.002.
[10] Pei B, Xu H, Xue Y. Lyapunov based estimation of flight stability boundary under icing conditions. Math Probl Eng 2017;2017:1–10. https://doi.org/10.1155/2017/6901894.
[11] Peyrl H, Parrilo PA. A theorem of the alternative for SOS Lyapunov functions. 2007 46th IEEE Conf. Decis. Control, IEEE; 2007, p. 1687–92. https://doi.org/10.1109/CDC.2007.4434258.
[12] Chesi G. Domain of attraction: analysis and control via SOS programming. vol. 415. London: Springer London; 2011. https://doi.org/10.1007/978-0-85729-959-8.
[13] Meng F, Wang D, Yang P, Xie G, Guo F. Application of sum-of-squares method in estimation of region of attraction for nonlinear polynomial systems. IEEE Access 2020;8:14234–43. https://doi.org/10.1109/ACCESS.2020.2966566.
[14] Kellett CM. Classical converse theorems in Lyapunov's second method. Discret Contin Dyn Syst - Ser B 2015;20:2333–60. https://doi.org/10.3934/dcdsb.2015.20.2333.
[15] Izumi S, Somekawa H, Xin X, Yamasaki T. Estimation of regions of attraction of power systems by using sum of squares programming. Electr Eng 2018;100:2205–16. https://doi.org/10.1007/s00202-018-0690-z.
[16] Li D, Ignatyev D, Tsourdos A, Wang Z. Region of attraction analysis for adaptive control of wing rock system. IFAC-PapersOnLine 2021;54:518–23. https://doi.org/10.1016/j.ifacol.2021.10.407.
[17] Li D, Ignatyev D, Tsourdos A, Wang Z. Nonlinear analysis for wing rock system with adaptive control. AIAA Scitech 2022 Forum, San Diego: AIAA; 2022, p. 8.
[18] Topcu U, Packard A. Local stability analysis for uncertain nonlinear systems. IEEE Trans Automat Contr 2009;54:1042–7. https://doi.org/10.1109/TAC.2009.2017157.

[19]   Genesio R, Tartaglia M, Vicino A. On the estimation of asymptotic stability regions: state of the art and new proposals. IEEE Trans Automat Contr 1985;30:747–55. https://doi.org/10.1109/TAC.1985.1104057.
[20]   Parrilo PA. Structured semidefinite programs and semialgebraic geometry methods in robustness and optimization. California Institute of Technology, Pasadena, CA, 2000. https://doi.org/10.7907/2K6Y-CH43.
[21]   Tan W, Packard A. Stability region analysis using sum of squares programming. 2006 Am. Control Conf., Minneapolis, MN, USA: IEEE; 2006, p. 6. https://doi.org/10.1109/ACC.2006.1656562.
[22]   Tibken B, Fan Y. Computing the domain of attraction for polynomial systems via BMI optimization method. Proc Am Control Conf 2006;2006:117–22. https://doi.org/10.1109/acc.2006.1655340.
[23]   Jarvis-Wloszek Z, Feeley R, Tan W, Sun K, Packard A. Some controls applications of sum of squares programming. 42nd IEEE Int. Conf. Decis. Control (IEEE Cat. No.03CH37475), vol. 5, IEEE; 2003, p. 4676–81. https://doi.org/10.1109/CDC.2003.1272309.
[24]   Kant N, Mukherjee R, Chowdhury D, Khalil HK. Estimation of the region of attraction of underactuated systems and its enlargement using impulsive inputs. IEEE Trans Robot 2019;35:618–32. https://doi.org/10.1109/TRO.2019.2893599.
[25]   Hafstein S, Giesl P. Review on computational methods for Lyapunov functions. Discret Contin Dyn Syst - Ser B 2015;20:2291–331. https://doi.org/10.3934/dcdsb.2015.20.2291.
[26]   Papachristodoulou A, Anderson J, Valmorbida G, Prajna S, Seiler P, Parrilo PA. Sum of squares optimization toolbox for Matlab user's guide 2016.
[27]   Prajna S, Papachristodoulou A, Parrilo PA. Introducing SOSTOOLs: a general purpose sum of squares programming solver. Proc. 41st IEEE Conf. Decis. Control. 2002., vol. 1, IEEE; 2002, p. 741–6. https://doi.org/10.1109/CDC.2002.1184594.
[28]   Seiler P. SOSOPT: a toolbox for polynomial optimization 2013;0:1–11.
[29]   Sturm JF. Using SeDuMi 1.02, a Matlab toolbox for optimization over symmetric cones. Optim Methods Softw 1999;11:625–53. https://doi.org/10.1080/10556789908805766.
[30]   Lofberg J. Pre- and post-processing sum-of-squares programs in practice. IEEE Trans Automat Contr 2009;54:1007–11. https://doi.org/10.1109/TAC.2009.2017144.
[31]   Tan W. Nonlinear control analysis and synthesis using sum-of-squares programming. Berkeley: University of California, 2006.
[32]   Bai Y, Wang Y, Svinin M, Magid E, Sun R. Function approximation technique based immersion and invariance control for unknown nonlinear systems. IEEE Control Syst Lett 2020;4:934–9. https://doi.org/10.1109/LCSYS.2020.2997600.
[33]   Prajna S, Papachristodoulou A, Seiler P, Parrilo PA. SOSTOOLS and its control applications. Lect. Notes Control Inf. Sci., vol. 312, 2005, p. 273–92. https://doi.org/10.1007/10997703_14.
[34]   Jarvis-Wloszek Z. Lyapunov based analysis and controller synthesis for polynomial systems using sum-of-squares optimization. Berkeley: University of California, 2003.
[35]   Chakraborty A, Seiler P, Balas GJ. Susceptibility of F/A-18 flight controllers to the falling-leaf mode: linear analysis. J Guid Control Dyn 2011;34:57–72. https://doi.org/10.2514/1.50674.
[36]   Tan W, Packard A. Searching for control Lyapunov functions using sums of squares programming. Allert Conf Commun Control Comput 2004:210–9.
[37]   Tan W, Packard A. Stability region analysis using polynomial and composite polynomial Lyapunov functions and sum-of-squares programming. IEEE Trans Automat Contr 2008;53:565–71. https://doi.org/10.1109/TAC.2007.914221.
[38]   Chesi G. Estimating the domain of attraction via union of continuous families of Lyapunov estimates. Syst Control Lett 2007;56:326–33. https://doi.org/10.1016/j.sysconle.2006.10.012.
[39]   Balestrino A, Crisostomi E, Grammatico S, Superiore S, Anna S, Martiri P, et al. R-composition of Lyapunov functions. 2009 17th Mediterr. Conf. Control Autom., IEEE; 2009, p. 126–31. https://doi.org/10.1109/MED.2009.5164527.

[40]   Balestrino A, Caiti A, Crisostomi E. Logical composition of Lyapunov functions. Int J Control 2011;84:563–73. https://doi.org/10.1080/00207179.2011.562549.
[41]   Chakraborty A, Seiler P, Balas GJ. Susceptibility of f/a-18 flight controllers to the falling-leaf mode: nonlinear analysis. J Guid Control Dyn 2011;34:73–85. https://doi.org/10.2514/1.50675.
[42]   Pandita R, Chakraborty A, Seiler P, Balas G. Reachability and region of attraction analysis applied to GTM dynamic flight envelope assessment. AIAA Guid. Navig. Control Conf., Reston, Virigina: AIAA; 2009. https://doi.org/10.2514/6.2009-6258.
[43]   Sidoryuk ME, Khrabrov AN. Estimation of regions of attraction of aircraft spin modes. J Aircr 2019;56:205–16. https://doi.org/10.2514/1.C034936.
[44]   Chakraborty A, Seiler P, Balas GJ. Nonlinear region of attraction analysis for flight control verification and validation. Control Eng Pract 2011;19:335–45. https://doi.org/10.1016/j.conengprac.2010.12.001.
[45]   Jones M, Mohammadi H, Peet MM. Estimating the region of attraction using polynomial optimization: a converse Lyapunov result. 2017 IEEE 56th Annu. Conf. Decis. Control, vol. 2018-Janua, IEEE; 2017, p. 1796–802. https://doi.org/10.1109/CDC.2017.8263908.
[46]   Shawki N, Lazarou G, Isenberg DR. Stability and performance analysis of a payload-manipulating robot without adaptive control. Int J Robot Autom 2020;35:23–34. https://doi.org/10.2316/J.2020.206-0077.
[47]   Dorobantu A, Seiler P, Balas G. Nonlinear analysis of adaptive flight control laws. AIAA Guid. Navig. Control Conf., Toronto, Ontario, Canada: AIAA; 2010. https://doi.org/10.2514/6.2010-8043.
[48]   Ryali V, Moudaglya KM. Robustness analysis of uncertain, nonlinear systems. Proc Am Control Conf 2000;5:3106–10. https://doi.org/10.1109/acc.2000.879137.
[49]   Davison EJ, Kurak EM. Computational method for determining quadratic Lyapunov functions for nonlinear systems. Automatica 1970;7:627–36. https://doi.org/10.1016/0005-1098(71)90027-6.
[50]   Ratschan S, She Z. Providing a basin of attraction to a target region of polynomial systems by computation of Lyapunov-like functions. SIAM J Control Optim 2010;48:4377–94. https://doi.org/10.1137/090749955.
[51]   Valmorbida G, Anderson J. Region of attraction estimation using invariant sets and rational Lyapunov functions. Automatica 2017;75:37–45. https://doi.org/10.1016/j.automatica.2016.09.003.
[52]   Wang S, She Z, Ge SS. Inner-estimating domains of attraction for nonpolynomial systems with polynomial differential inclusions. IEEE Trans Cybern 2020:1–14. https://doi.org/10.1109/TCYB.2020.2987326.

# Figures

(Color should be used in print please.)

Fig. 1 Demonstration of shape functions: $p_1(x): N=[1,0;0,1], x^*=[0,0]$ ; $p_2(x): N=[1,1;0,3], x^*=[0,0]$ ; $p_3(x): N=[1,0;0,3], x^*=[1,1]$

Fig. 2 Conceptual diagram of three rounds of shifts in RcomSSF

Fig. 3 Location of shifting centers

Fig. 4 Two ellipses(blue lines); Left: the union between two ellipses $\hat{R}_1 \cup \hat{R}_2 = \{x \in \mathbb{R}^n : R_\cup(R_1, R_2) > 0\}$ (green line); Right: the intersection between two ellipses $\hat{R}_1 \cap \hat{R}_2 = \{x \in \mathbb{R}^n : R_\cap(R_1, R_2) > 0\}$ (red line)

Fig. 5 ROA estimation for Example 1

  (a) ROA estimation by A1 and A2 at 30 and 60 iterations

  (b) ROA estimation and level sets of $p(x)$ by A2

  (c) Shifting details of RcomSSF

  (d) ROA estimation by RcomSSF

Fig. 6 ROA estimation for Example 2

  (a) ROA estimation by A1 and A2

  (b) ROA estimation by RcomSSF

Fig. 7 ROA estimation for Example 3

  (a) Vector field and ROA estimation by A1, A2, and RcomSSF

  (b) Effect of shifting centers

  (c) Effect of shape matrices

Fig. 8 ROA estimation for Example 4

(a) ROA estimation by A1, A2, RcomSSF, and methods in [28,43]

(b) Shifting details of RcomSSF

Fig. 9 ROA estimation for Example 5

(a) ROA estimation by A1 and A2

(b) ROA estimation by RcomSSF

(c) Crosssection for $x_2 = 0$

(d) Shifting details of RcomSSF

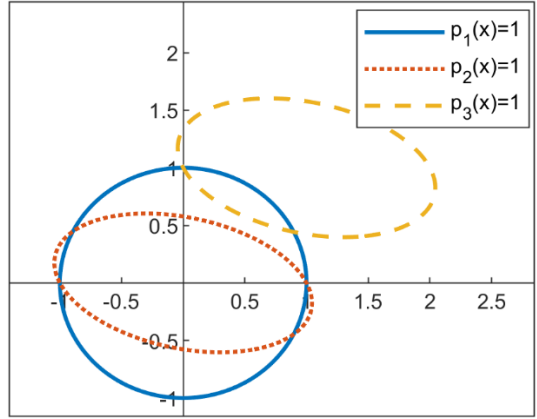

**Fig. 1 Demonstration of shape functions:** $p_1(x): N = [1,0;0,1], x^* = [0,0]$; $p_2(x): N = [1,1;0,3], x^* = [0,0]$; $p_3(x): N = [1,0;0,3], x^* = [1,1]$

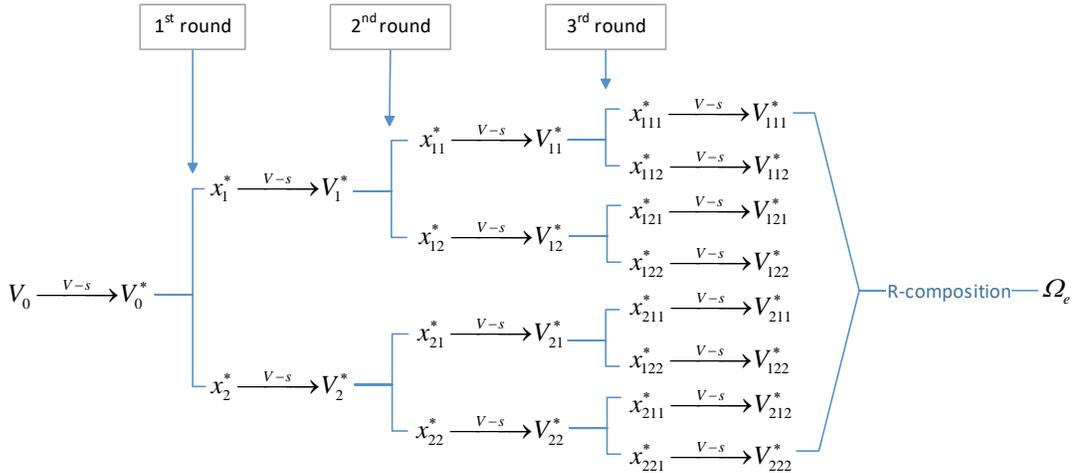

**Fig. 2 Conceptual diagram of three rounds of shifts in RcomSSF**

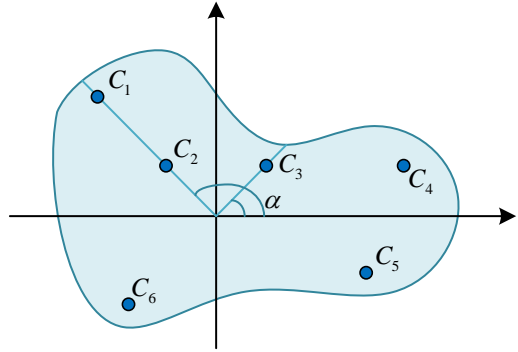

**Fig. 3 Location of shifting centers**

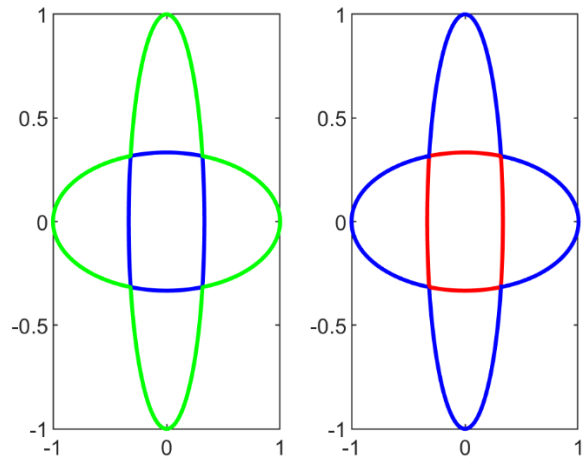

**Fig. 4 Two ellipses(blue lines);**
**Left: the union between two ellipses $\hat{R}_1 \bigcup \hat{R}_2 = \{x \in \mathbb{R}^n : R_{\bigcup}(R_1, R_2) > 0\}$ (green line);**
**Right: the intersection between two ellipses $\hat{R}_1 \bigcap \hat{R}_2 = \{x \in \mathbb{R}^n : R_{\bigcap}(R_1, R_2) > 0\}$ (red line).**

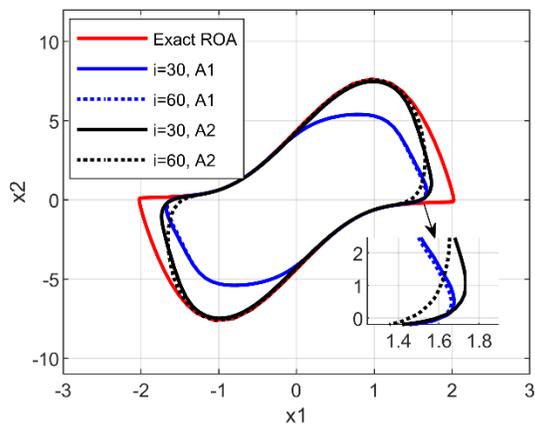

(a) ROA estimation by A1 and A2 at 30 and 60 iterations

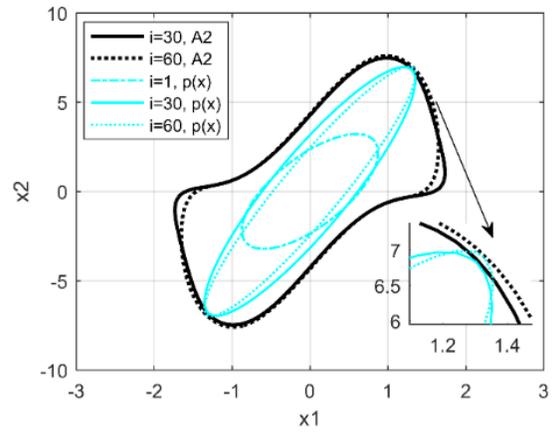

(b) ROA estimation and level sets of $p(x)$ by A2

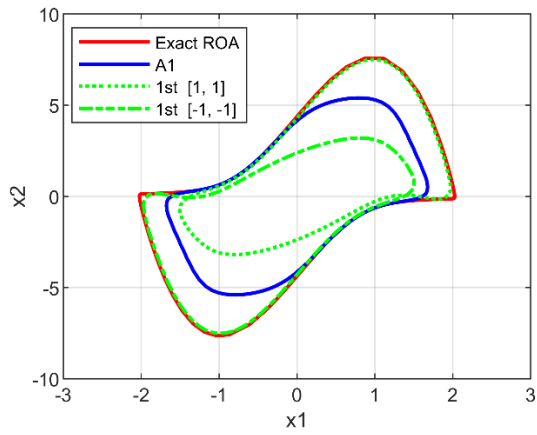

(c) Shifting details of RcomSSF

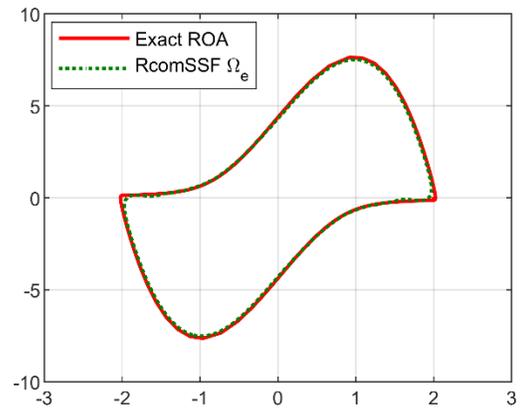

(d) ROA estimation by RcomSSF

**Fig. 5 ROA estimation for Example 1**

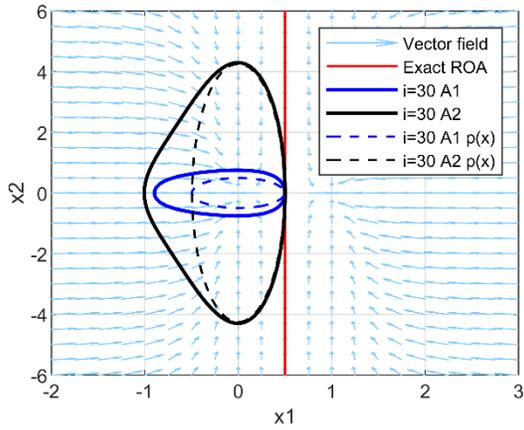
(a) ROA estimation by A1 and A2

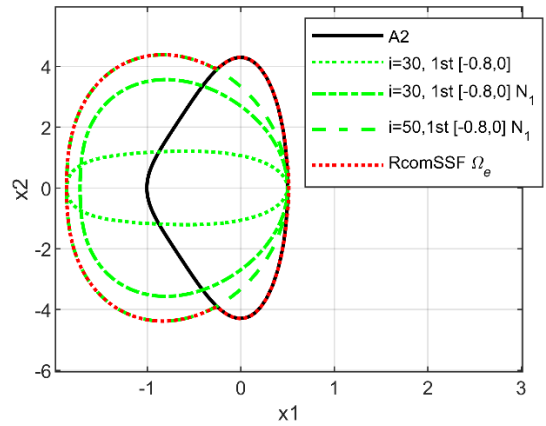
(b) ROA estimation by RcomSSF

**Fig. 6 ROA estimation for Example 2**

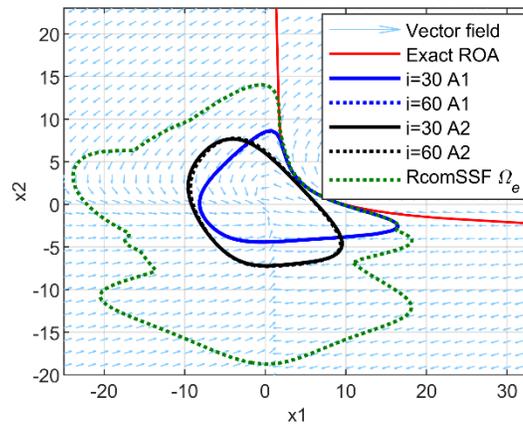
(a) Vector field and ROA estimation by A1, A2, and RcomSSF

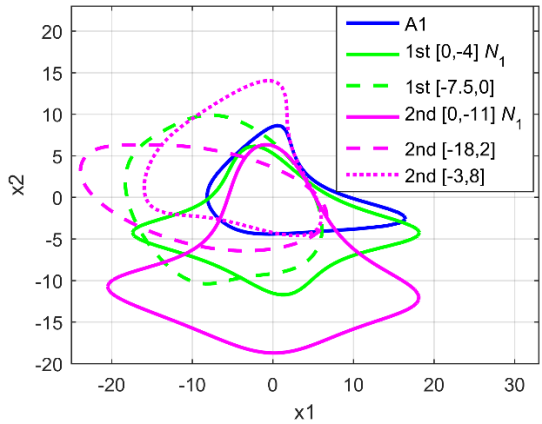
(b) Effect of shifting centers

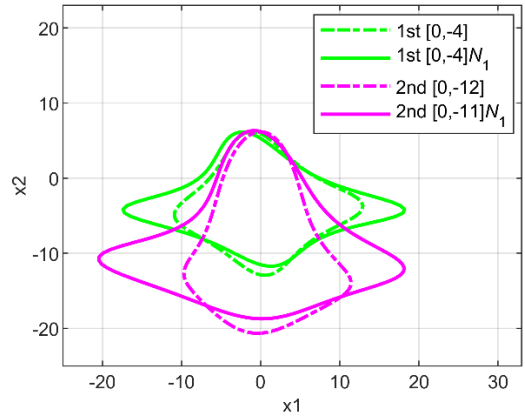
(c) Effect of shape matrices

**Fig. 7 ROA estimation for Example 3**

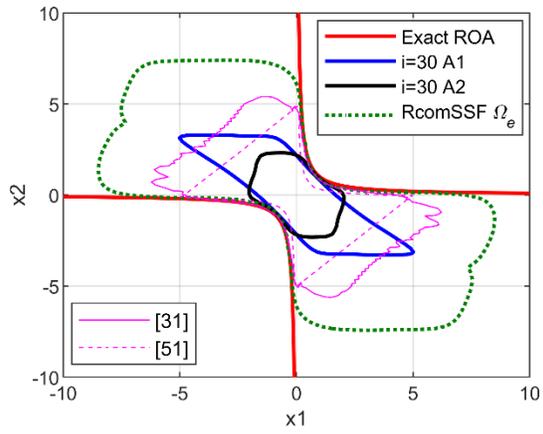 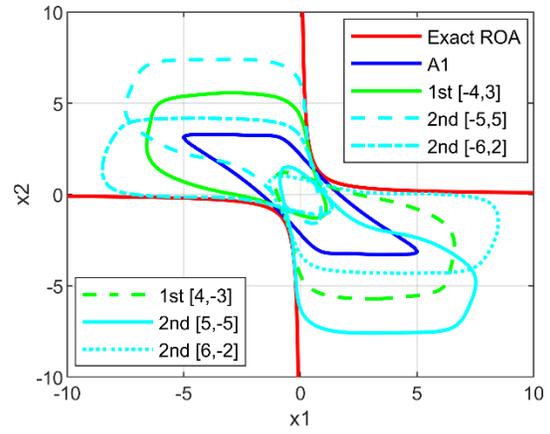

(a) ROA estimation by A1, A2, RcomSSF, and methods in [31,51]

(b) Shifting details of RcomSSF

**Fig. 8 ROA estimation of Example 4**

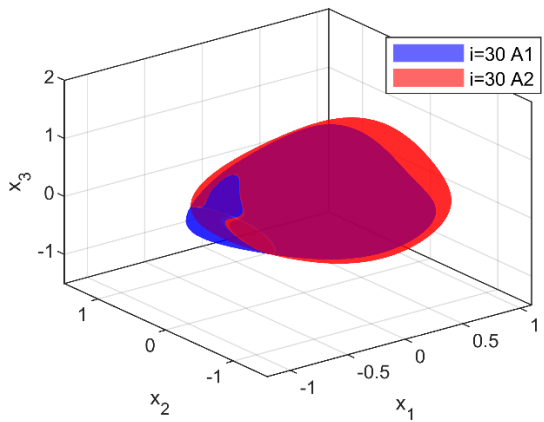 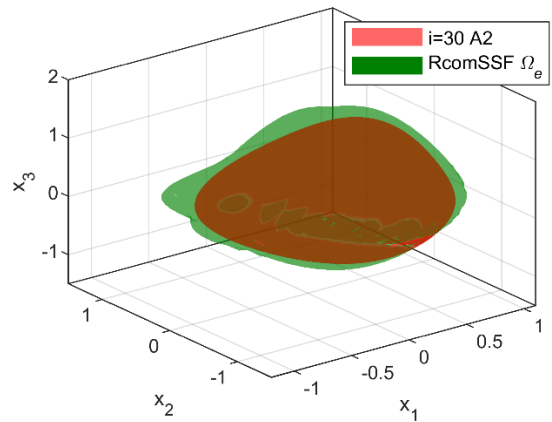

(a) ROA estimation by A1 and A2

(b) ROA estimation by RcomSSF

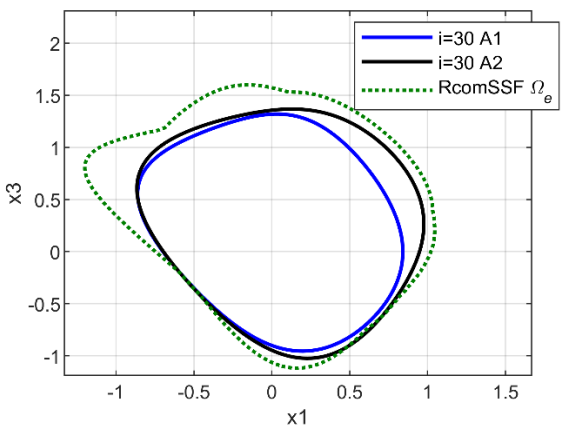 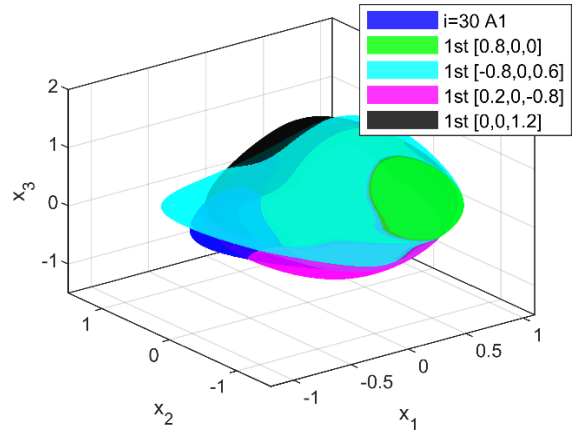

(c) Crosssection for $x_2 = 0$

(d) Shifting details of RcomSSF

**Fig. 9 ROA estimation for Example 5**

# Tables

Table 1  Correspondence between logical functions and R-composition

| Boolean | Geometry | R-composition |
|---|---|---|
| **not** | complement | $-R(x)$ |
| **and** | intersection | $R_\cap(R_1, R_2) = R_1(x) + R_2(x) - \sqrt{R_1^2(x) + R_2^2(x) - \tau R_1(x) R_2(x)}$ |
| **or** | union | $R_\cup(R_1, R_2) = R_1(x) + R_2(x) + \sqrt{R_1^2(x) + R_2^2(x) - \tau R_1(x) R_2(x)}$ |